\documentclass[12pt]{article}

\usepackage{mathtools}
\usepackage{amsmath,amssymb,amsbsy}
\usepackage{graphicx}
\usepackage{multirow}

\newtheorem{thm}{Theorem}[section]
\newtheorem{lem}[thm]{Lemma}
\newtheorem{obs}[thm]{Observation}

\newtheorem{prob}[thm]{Problem}
\newtheorem{dfn}[thm]{Definition}

\newtheorem{conj}[thm]{Conjecture}

\makeatletter
\def\imod#1{\allowbreak\mkern10mu({\operator@font mod}\,\,#1)}
\makeatother

\newcommand{\wlg}{we may suppose without loss of generality }

\newcommand{\zet}{\mathbb{Z}}
\def\gr{\mathcal{G}}

\newcommand{\qed}{\hfill \rule{.1in}{.1in}}

\def\gcd{\mathop{\rm gcd}\nolimits}

\begin{document}
\title{A magic rectangle  set on  Abelian groups}

\author{Sylwia Cichacz$^{1,3}$,  Tomasz Hinc$^{2,3}$\\
$^1$AGH University of Science and Technology Krak\'ow, Poland\\
$^2$Pozna\'n University of Economics and Business, Poland\\
$^3$University of Primorska, Slovenia}

\maketitle
\begin{abstract}

A $\Gamma$-magic rectangle  set $MRS_{\Gamma}(a, b; c)$  of order $abc$ is a collection of $c$ arrays $(a\times b)$
whose entries are elements of group $\Gamma$, each appearing once, with all row sums in
every rectangle equal to a constant $\omega\in \Gamma$ and all column sums in every rectangle equal to a
constant $\delta \in \Gamma$. \\

In this paper we prove that for  $\{a,b\}\neq\{2^{\alpha},2k+1\}$ where $\alpha$ and $k$ are some natural numbers, a  $\Gamma$-magic rectangle set MRS$_{\Gamma}(a, b;c)$ exists if and only if $a$ and $b$ are both even or and $|\Gamma|$ is odd or   $\Gamma$ has more than one involution. Moreover we obtain  sufficient and necessary conditions for existence a $\Gamma$-magic rectangle MRS$_{\Gamma}(a, b)$=MRS$_{\Gamma}(a, b;1)$.

\end{abstract}

 \section{Introduction}
 A \textit{magic square} of order $n$ is an $n\times n$ array with entries $1,2,\ldots,n^2$, each
appearing once, such that the sum of each row, column, and both main diagonals
is equal to $n(n^2+1)/2$.  The earliest known magic square is a $3\times3$ magic square called \textit{Lo Shu magic square} can be traced in Chinese literature as far back as 2800 B.C. Since then, certainly, many people studied magic squares. For a  survey of magic squares, see Chapter 34 in \cite{ColDin}. Magic rectangles are a natural generalization of magic squares. 
A {\it magic rectangle} MR$(a, b)$ is an $a \times b$ array with entries from the set
$\{1, 2, \ldots , ab\}$, each appearing once, with all its row sums equal to a constant $\delta$ and with
all its column sums equal to a constant $\eta$. It was proved in \cite{Har1,Har2}:

\begin{thm}[\cite{Har1,Har2}] A magic rectangle MR$(a, b)$ exists if and only if $a, b > 1$, $ab > 4$,
and $a \equiv b \pmod 2$.\label{Har}
\end{thm}
The  generalization of magic rectangles  was introduced by Froncek in \cite{Fro}.  Namely, a \textsl{magic rectangle set} M= MRS$(a, b; c)$ is a collection of $c$ arrays ($a\times b$)
whose entries are elements of $\{1, 2, \ldots , abc\}$, each appearing once, with all row sums in
every rectangle equal to a constant $\delta$ and all column sums in every rectangle equal to a
constant $\eta$.

It was shown the following.
\begin{thm}[\cite{Fro2}] For $a > 1$ and $b \geq 4$, a magic rectangle set MRS$(a, b; c)$ exists
if and only if $a, b \equiv 0 \pmod 2$ or $abc \equiv 1 \pmod 2$.\label{thmFro}
\end{thm}

Assume $\Gamma$ is an Abelian group of order $n$ with the operation denoted by $+$.  For convenience
we will write $ka$ to denote $a + a + \ldots + a$ (where the element $a$ appears $k$ times), $-a$ to denote the inverse of $a$ and
we will use $a - b$ instead of $a+(-b)$.  Moreover, the notation $\sum_{a\in S}{a}$ will be used as a short form for $a_1+a_2+a_3+\dots$, where $a_1, a_2, a_3, \dots$ are all elements of the set $S$. The identity element of $\Gamma$ will be denoted by $0$. Recall that any group element $\iota\in\Gamma$ of order 2 (i.e., $\iota\neq 0$ and $2\iota=0$) is called an \emph{involution}.

Evans introduced modular magic rectangles \cite{Eva}. A \textit{modular rectangle}   is an  $(a\times b)$ array with entries $1,2,\ldots,ab$ such that  all the rows sum are congruent mudulo $ab$ and all the rows sum are congruent modulo $ab$.  The concept of magic squares on general Abelian groups was presented in~\cite{SY}. A {\it $\Gamma$-magic square} MS$_{\Gamma}(n)$ is an $n \times n$ array with entries from the an Abalian group $\Gamma$ of order $n^2$, each appearing once, with all its row, column and diagonal sums equal to a constant $\delta$. It was proved in \cite{SY}:

\begin{thm}[\cite{SY}] $\Gamma$-magic squares MS$_{\Gamma}(n)$ exist for all groups $\Gamma$ of order $n^2$ for any $n > 2$.\label{SY}
\end{thm}

The author has introduced  the following generalization of $\Gamma$-magic squares MS$_{\Gamma}(n)$ (\cite{ref_CicIWOCA}).
\begin{dfn}
A $\Gamma$-\textit{magic rectangle  set} MRS$_{\Gamma}(a, b; c)$ on group $\Gamma$ of order $abc$ is a collection of $c$ arrays $(a\times b)$
whose entries are elements of group $\Gamma$, each appearing once, with all row sums in
every rectangle equal to a constant $\omega\in \Gamma$ and all column sums in every rectangle equal to a
constant $\delta \in \Gamma$. If $c=1$, then  a $\Gamma$-magic rectangle  set MRS$_{\Gamma}(a, b; 1)$ we call a $\Gamma$-\textit{magic rectangle}   MR$_{\Gamma}(a, b)$ 
\end{dfn}
It was proved the following.
\begin{thm}[\cite{ref_CicIWOCA}]
If $a,b$ are both even, then a magic rectangle set MRS$_\Gamma(a, b; c)$ exists
for every $c$ for any $\Gamma$ of order $abc$.\label{glowneIWOCA}
\end{thm}

\begin{obs}[\cite{ref_CicIWOCA}]\label{codd}  If $a$ is even, $b$ is odd then for any $c$ and an Abelian group $\Gamma$ having exactly one involution, $|\Gamma|=abc$ there does not exist a $\Gamma$-magic rectangle set MRS$_{\Gamma}(a, b; c)$.
\end{obs}
\begin{obs}[\cite{ref_CicIWOCA}]\label{odd} If exactly one of numbers the $a,b,c$ is even, then  there does not exist a magic rectangle set MRS$_{\Gamma}(a, b; c)$ on any Abelian group $\Gamma$ having exactly one involution.
\end{obs}


The paper is organized as follows.  In Section 2 we give some preliminaries on Abelian groups. In Section~3  we show that  for  $\{a,b\}\neq\{2^{\alpha},2k+1\}$ where $\alpha$ and $k$ are some natural numbers, a  $\Gamma$-magic rectangle set MRS$_{\Gamma}(a, b;c)$ exists if and only if $a$ and $b$ are both even or $\Gamma\in\gr$. This result implies  that if a MRS$(a, b; c)$ exists then a MRS$_{\Gamma}(a, b; c)$ exists for any Abelian group $\Gamma$ of order $abc$. Furthermore we apply this result for giving sufficient and necessary conditions for existence a $\Gamma$-magic rectangle MR$_{\Gamma}(a, b)$.  In the last section we post some open problems.
\section{Preliminaries}

 A non-trivial
finite group has elements of order $2$ (involutions) if and only if the order of the group is even. The fundamental theorem of finite Abelian groups states that a finite Abelian
group $\Gamma$ of order $n$ can be expressed as the direct product of cyclic subgroups of prime-power order. This implies that
$$\Gamma\cong\zet_{p_1^{\alpha_1}}\times\zet_{p_2^{\alpha_2}}\times\ldots\times\zet_{p_k^{\alpha_k}}\;\;\; \mathrm{where}\;\;\; n = p_1^{\alpha_1}\cdot p_2^{\alpha_2}\cdot\ldots\cdot p_k^{\alpha_k}$$
and $p_i$ for $i \in \{1, 2,\ldots,k\}$ are not necessarily distinct primes. This product is unique up to the order of the direct product. When $p$ is the number of these cyclic components whose order is a multiple of $2$, then $\Gamma$ has $2^p-1$ involutions. In particular, if $n \equiv 2 \pmod 4$, then $\Gamma\cong \zet_2\times \Lambda$ for some
Abelian group $\Lambda$ of odd order $n/2$.  Moreover every cyclic group of an even order has exactly one involution. \\

One of the consequences of the fundamental theorem of finite Abelian groups is that  for any divisor $m$ of $|\Gamma|$ there exists a subgroup $H$ of $\Gamma$ of order $m$.  A \textit{quotient group} $\Gamma$ \textit{modulo} $H$ for a subgroup $H$ of $\Gamma$ will be denoted by $\Gamma/H$. It is wellknown that $|\Gamma/H|=|\Gamma|/|H|$. \\

For convenience, let $\gr$ denote the set consisting of all Abelian groups which are
of odd order or contain more than one involution. Since the properties and results in this paper are invariant under the isomorphism between groups, we only need to consider one group in each isomorphism class.

\section{A $\Gamma$-magic rectangle  set MRS$_{\Gamma}(a, b; c)$}
In \cite{ref_Wal} Marr and Wallis gave a definition of a Kotzig array as a $j\times k$ grid, each row being a permutation of $\{0,1,\dots,k-1\}$ and each column having the same sum. Note that a  Latin square  is a special case of a Kotzig array. Moreover a Kotzig array 
has the “magic-like” property with respect to the sums of rows and columns, but at the same time allows entry repetitions in columns. Kotzig arrays play an important role for graph labelings.  There are many constructions based on Kotzig arrays of various magic-type constructions \cite{ref_CicFroSin, ref_CicNik,Fro2,ref_InaLlaMor,ref_McQ,ref_Wal}.

  In the paper \cite{ref_CicZ} the author introduced a generalization of Kotzig arrays and gave necessary and sufficient conditions for their existence. Namely, for an Abelian group $\Gamma$ of order $k$ we define a $\Gamma$-\textit{Kotzig array}  of size $j\times k$ as a  $j\times k$ grid, where each row is a permutation of elements of $\Gamma$ and each column has the same sum. Moreover \wlg that all the elements in the first column are  $0$. 

\begin{thm}[\cite{ref_CicZ}]\label{Kotzig}
A $\Gamma$-Kotzig array of size $j \times k$ exists if and only if
$j>1$ and $j$ is even or $\Gamma\in \gr$. 
\end{thm}
We we start with a lemma which plays an important role in the proof of the
main result.

\begin{lem}\label{lemgl}
Let $\Gamma$ be a group of order $abc$. If there exists a subgroup $\Gamma_0$ of order $ab$ such that $\Gamma/\Gamma_0$ belongs to $\gr$ and there exists a $\Gamma_0$-magic rectangle  MR$_{\Gamma_0}(a, b)$, then there exists a $\Gamma$-magic rectangle set MRS$_{\Gamma}(a, b; c)$.
\end{lem}
\textit{Proof.} If $a,b$ are both even, then a $\Gamma$-magic rectangle MR$_\Gamma(a, b)$ exists by Theorem~\ref{glowneIWOCA}. Observe that if  MR$_{\Gamma}(a, b)$ exists, then obviously  MR$_{\Gamma}(b, a)$ exists as well thus \wlg that $a$ is odd. Moreover if $b$ is also odd, then  without loss of generality we assume that $a\leq b$. \\
Denote  by $x_{i,j}$ the $i$-th row and $j$-th column of the $\Gamma_0$-magic rectangle MR$_{\Gamma_0}(a, b)$ and  $\omega_0$ and $\delta_o$ are the row sum and the column sum, respectively in MR$_{\Gamma_0}(a,b)$.\\
Note that $|\Gamma/\Gamma_0|=c$. Based on Theorem~\ref{Kotzig} there exists $\Gamma/\Gamma_0$-Kotzig array of size $a\times c$, where the entries are the coset representatives (and one of them is $0$) for $\Gamma/\Gamma_0$.  Denote by $k_{i,s}$ the $i$-th row and $s$-th column of the $\Gamma/\Gamma_0$-Kotzig array. Let all column sums be $\sum_{i=1}^ak_{i,s}=0$. 
Using Kotzig array we build $c$ different $a\times b$ \textit{residual rectangles} $R^s(a, b)$ with entries $r^s_{i,j}$
for $1 \leq s \leq c$. In the first column of a given $R^s(a, b)$, we place the $s$-th column of $\Gamma/\Gamma_0$-Kotzig array. That is, $r^s_{i,1} =k_{i,s}$. We will consider two cases on the parity of $b$.\\

\textit{Case 1.} $b$ is odd.\\
In the following $a-1$ columns, we place a circulant array constructed from the first column.  If $b>a$, then all remaining even columns will be the same as the first column, while the remaining odd columns will be filled with the complements
of the entries in the previous column with respect to $\Gamma/\Gamma_0$. Namely, set 
$$r^s_{i,j}=\left\{
\begin{array}{lcl}
r^s_{(i+j-1)(\mathrm{mod}\, a),1}, & \text{if} & j =2,3,\ldots,a, \\
r^s_{i,1}, & \text{if} & j =a+1,a+3,\ldots,b-1, \\
-r^s_{i,1}, & \text{if} & j =a+2,a+4,\ldots,b.
\end{array}%
\right.
$$

\textit{Case 2.} $b$ is even.\\
In this case all odd columns will be the same as the first column, whereas the even columns will be filled with the complements
of the entries in the previous column with respect to $\Gamma/\Gamma_0$. Namely, set 
$$r^s_{i,j}=\left\{
\begin{array}{lcl}
r^s_{i,1}, & \text{if} & j =3,5,\ldots,b-1, \\
-r^s_{i,1}, & \text{if} & j =2,4,\ldots,b.
\end{array}%
\right.
$$

\vspace{0.5 cm}

Observe that in both cases  $r^s_{i,j}\neq r^{s'}_{i,j}$ for any $s\neq s'$ and moreover  the sum of every  column is  $\sum_{j=1}^ar^s_{i,j}=\sum_{j=1}^ak_{s,j}=0$ and the sum of each row is 
 $\sum_{i=1}^br^s_{i,j}=0$. 

Denote by $x^s_{i,j}$ the entry in the $i$-th row and $j$-th column of the $s$-th rectangle in the $\Gamma$-magic rectangle set MRS$_{\Gamma}(a, b; c)$. Let $x_{i,j}^s=x_{i,j}+r^s_{i,j}$, for $i=1,2,\ldots,a$, $j=1,2,\ldots,b$ and $s=1,2,\ldots,k$.

Obviously  all  column sums equal to a
constant $\delta=\delta_0$ whereas row sums equal to $\omega=\omega_0$.  
We still need to verify that every  element of $\Gamma$ appears exactly once. Suppose that $x_{i,j}^s= x_{i',j'}^{s'}$. If $(i, j)\neq(i',j')$, then  $x_{i,j}\neq x_{i',j'}$ what implies that also  $r^s_{i,j}\neq r^{s'}_{i',j'}$. On one hand  $0\neq r^s_{i,j}- r^{s'}_{i',j'}\not\in\Gamma_{0}$, on the other $r^s_{i,j}- r^{s'}_{i,j}=x_{i,j}-x_{i',j'}\in\Gamma_{0}$, a contradiction. Assume now $(i, j)=(i',j')$, recall that by definition of residual rectangles we have  $r^s_{i,j}\neq r^{s'}_{i,j}$  and this implies that again  $0\neq r^s_{i,j}- r^{s'}_{i,j} =x_{i,j}-x_{i,j}=0$, a contradiction.~\qed

\begin{obs}\label{obs1}
Let $p_1,p_2\geq3$ be primes. A $\Gamma$-magic rectangle  MR$_{\Gamma}(p_1, p_2)$ exists on any group $\Gamma$ of order $p_1p_2$.
\end{obs}
\textit{Proof.} Assume first that  $\Gamma\cong \zet_{p_1p_2}$. Since  there exists MR$(p_1, p_2)$ by Theorem~\ref{Har}, therefore there exists a MR$_{\Gamma}(p_1, p_2)$. If now $\Gamma\not\cong \zet_{p_1p_2}$, then because $p_1$ and $p_2$ are primes note that $p_1=p_2$ and  $\Gamma\cong \zet_{p_1}\times \zet_{p_1}$. Denote by $x_{i,j}$ the entry in the $i$-th row and $j$-th column of the  rectangle. Set $x_{i,j}=(i,j)$. In this  case every column and every row adds up to $(0,0)=0$.~\qed

\begin{lem}\label{lem2}
 If there exist odd  primes $p_1,p_2$ such that $p_1|a$ and $p_2|b$ and moreover  $\Gamma\in \gr$ is a group of order $abc$, then there exists a $\Gamma$-magic rectangle set  MRS$_{\Gamma}(a, b;c)$.
\end{lem}
\textit{Proof.} 
Note that since $\Gamma$ is an Abelian group there exists a subgroup $\Gamma_0$  of $\Gamma$ of order $p_1p_2$.  By  Observation~\ref{obs1} there exists a $\Gamma_0$-magic rectangle  MR$_{\Gamma_0}(p_1, p_2)$. Moreover since $p_1,p_2\geq3$ are odd primes the quotient group $\Gamma/\Gamma_0$ belongs to $\gr$. Recall that $|\Gamma/\Gamma_0|=abc/p_1p_2$ and thus Lemma~\ref{lemgl} implies that there exists a  MRS$_{\Gamma}(p_1, p_2;abc/p_1p_2)$. To construct  one of  the rectangles from MRS$_{\Gamma}(a, b;c)$ we simply take  $ab/p_1p_2$ of  MRS$_{\Gamma}(p_1, p_2;abc/p_1p_2)$ rectangles and ''glue'' them  into a  rectangle $a \times b$.~\qed\\


Now we can state the main result.

\begin{thm}\label{main2} Let   $\{a,b\}\neq\{2^{\alpha},2l+1\}$. A  $\Gamma$-magic rectangle set MRS$_{\Gamma}(a, b;c)$ exists if and only if $a$ and $b$ are both even or $\Gamma\in\gr$.
\end{thm}
\textit{Proof.} If $a$ and $b$ are both even  then we are done by Theorem~\ref{glowneIWOCA}. 
Suppose  that $a$ and $b$ are both odd. For $\Gamma\in \gr$  we apply Lemma~\ref{lem2}. If  $\Gamma\not\in \gr$, then $\Gamma$ has exactly one involution and therefore $c$ has to be even. As a consequence there does not  a MRS$_{\Gamma}(a, b;c)$ by Observation~\ref{odd}.

Assume now  $a$ is odd and $b=(2m+1)2^{\alpha}$ for some integers $m>0$, $\alpha>0$. 
Since $|\Gamma|=abc\equiv 0 \pmod 2$ the group $\Gamma$ has to contain at least one involution. If $\Gamma$ has exactly one involution then by Observation~\ref{codd}  there does not exist a $\Gamma$-magic rectangle  set MRS$_{\Gamma}(a, b;c)$. Assume now that $\Gamma$ has more than one involution, thus $\Gamma \in \gr$. Let $p_1,p_2$ be odd primes such that $p_1|a$ and $p_2|b$, hence Lemma~\ref{lem2} implies that there exists a $\Gamma$-magic rectangle set MRS$_{\Gamma}(a, b;c)$.~\qed\\

The following theorem immediately follows from the 	above theorem.
\begin{thm}\label{main} Let $a > 1$ and $b \geq 4$.
If  a magic rectangle  set MRS$(a, b; c)$ exists then a $\Gamma$-magic rectangle  set  MRS$_{\Gamma}(a, b; c)$ exists for any Abelian group $\Gamma$ of order $abc$.
\end{thm}
\textit{Proof.} For $a > 1$ and $b \geq 4$, a magic rectangle set MRS$(a, b; c)$ exists
if and only if $a, b \equiv 0 \pmod 2$ or $abc \equiv 1 \pmod 2$ by Theorem~\ref{thmFro}.
If  $abc \equiv 1 \pmod 2$, then $\Gamma \in \gr$. Hence we  are done by Theorem~\ref{main2}.~\qed\\

Note that in and Theorem~\ref{main2} the only  missing case in the full characterization of $\Gamma$-magic rectangle sets  MRS$_{\Gamma}(a, b;c)$ is the case when $a$ is odd, $b=2^{\alpha}>1$. So far we are able to prove the following.

\begin{lem}\label{p22} 
Let $p$ be a prime. There exists a  $\Gamma$-magic rectangle   MR$_{\Gamma}(p, 2^{\alpha})$ if and only if a group $\Gamma$ has more than one involution.\end{lem}
\textit{Proof.}

By Theorem~\ref{glowneIWOCA} we can assume that $p\geq 3$. Moreover since $|\Gamma|$ is even the group $\Gamma$ has at least one involution. If it has exactly one involution then there does not exits a MRS$_{\Gamma}(p, 2^{\alpha},1)$=MR$_{\Gamma}(p, 2^{\alpha})$ by Observation~\ref{codd}. Therefore $\Gamma\cong\zet_p\times \zet_{2^{\alpha_1}}\times\zet_{2^{\alpha_2}}\times \ldots \times \zet_{2^{\alpha_k}}$ for $k>1$. Without loss of generality we can assume  $1\leq \alpha_1\leq \alpha_2\leq\ldots \leq\alpha_k$.  Let $\Delta=\zet_{2^{\alpha_1}}\times\zet_{2^{\alpha_2}}\times \ldots \times \zet_{2^{\alpha_k}}$ and $\alpha=\sum_{i=1}^k\alpha_i$.  We will consider two cases on $p$.\\

\textit{Case 1.} $\gcd(2^{\alpha}-1,p)=1$\\
Let $f(x)=-(2^{\alpha}-1)x$ for any $x\in \zet_p$. Note that since $\gcd(2^{\alpha}-1,p)=1$ the mapping $f$ is an automorphism. Denote by $k_{i,j}$ the $i$-th row and $j$-th column of the $\Delta$-Kotzig array of size of size $p\times 2^{\alpha}$ (which exists by Lemma~\ref{Kotzig}).\\

Denote  by $x_{i,j}$ the $i$-th row and $j$-th column of the $\Gamma$-magic rectangle MR$_{\Gamma}(p, 2^{\alpha})$.
Define $x_{i,j}=(a_{i,j},b_{i,j})$ for $i=1,2,\ldots,p$, $j=1,2,\ldots,2^{\alpha}$ such that $a_{i,j}\in \zet_p$ and $b_{i,j}\in \Delta$. Set
$$x_{i,j}=\left\{
\begin{array}{lcl}
(f(i-1),k_{i,j}), & \text{if} & j=1, \\
(i-1,k_{i,j}), & \text{if} & j\in\{2,3,\ldots,2^{\alpha}), \\
\end{array}%
\right.
$$
for $i=1,2,\ldots,p$.

In this  case every column sum is $(\sum_{x\in \zet_p}x,\sum_{i=1}^pk_{i,j})=(0,0)$ and every row sum  is $(0,0)$.\\

\textit{Case 2.} $\gcd(2^{\alpha}-1,p)=p$\\
We will create now a  $\Delta$-Kotzig array of size of size $p\times 3$ with the entry in the $i$-th row and $j$-th column denoted by $k_{i,j}$. Let $k'_{1,j}=k_{1,j}-k_{1,2}+k_{1,3}$ and $k'_{i,j}=k_{i,j}$ for $i=2,3$. Let the  first $3$ rows in MR$_{\Gamma}(p, 2^{\alpha})$ be given in the table:

\begin{center}
\begin{tabular}{|c|c|c|c|c|c|c|}
  \hline
	$(1, k_{1,1}')$& $(0, k_{3,2}')$  & $(1, k_{1,3}')$& $(1, k_{1,4}')$&$(1, k_{1,5}')$&\ldots&$(1, k_{1,2^{\alpha}}')$\\ \hline
	$(0, k_{3,1}')$& $(p-1, k_{2,2}')$& $(p-1, k_{2,3}')$& $(p-1, k_{2,4}')$&$(p-1, k_{2,5}')$&\ldots&$(p-1, 2^{\alpha})$\\ \hline
	$(p-1,k_{2,1}')$& $(1,k_{1,2}')$  & $(0, k_{3,3}')$& $(0, k_{3,4}')$& $(0, k_{3,5}')$&\ldots& $(0, 2^{\alpha})$\\ \hline
\end{tabular}
\end{center}
All other than $0,1,p-1$ elements from $\zet_p$  we join into pairs $(x,-x)$ (note that because $p$ is odd we have $x\neq -x$) and create two rows MR$_{\Gamma}(p, 2^{\alpha})$ as follows:
\begin{center}
\begin{tabular}{|c|c|c|c|c|c|c|c|}
  \hline
$(x, k_{1,1}')$& $(-x, k_{1,2}')$& $(x,k_{1,3}')$&$(-x,k_{1,4}')$&\ldots&$(x,k_{1,2^{\alpha}-1}')$&$(-x,k_{1,2^{\alpha}}')$\\\hline
$(-x, k_{1,1}')$& $(x, k_{1,2}')$&$(-x,k_{1,3}')$&$(x,k_{1,4}')$&\ldots&$(-x,k_{1,2^{\alpha}-1}')$&$(x,k_{1,2^{\alpha}}')$\\\hline
\end{tabular}
\end{center}

Note that  $k'_{2,1}=k'_{3,1}$ (since $k_{1,1}=k_{2,1}=k_{3,1}=0$) and moreover $k'_{1,2}=k'_{3,2}$ and this implies that every  element of $\Gamma$ appears exactly once in the table. Moreover one can easily see that  every column sum is $(\sum_{x\in \zet_p}x,\sum_{i=1}^pk_{i,j})=(0,-k_{1,2}+k_{1,3})$ and every row sum  is $(0,0)$. This finishes the proof.~\qed\\

\begin{thm}\label{a22} Let $a$ be odd, $b=2^{\alpha}>2$ and  $c=2^{\beta}(2m+1)$. Let $\Gamma\cong \Delta \times \zet_{2^{\alpha_1}}\times\zet_{2^{\alpha_2}}\times \ldots \times \zet_{2^{\alpha_k}}$ for  $\sum_{i=1}^k\alpha_i=\alpha+\beta$, $k>1$ and some group $\Delta$ of odd order $a(2m+1)$. Assume that  $1\leq \alpha_1\leq \alpha_2\leq\ldots \leq\alpha_k$. There exists  a MRS$_{\Gamma}(a, 2^{\alpha};c)$ if one of the following conditions hold:
\begin{enumerate}
	\item $\alpha_1>1$,
	\item $\alpha_3>1$,
	\item $k>3$,
		\item $c$ is odd.
\end{enumerate}
\end{thm}
\textit{Proof.} Note that  $b\equiv 0 \pmod 4$. Let $p$ be an odd prime such that $p|a$. \\

Assume now that $\alpha_1>1$  or $\alpha_3>1$ or $k\geq4$. In all three cases there exists a subgroup $\Gamma_0$ of $\Gamma$ such that $\Gamma_0\cong \zet_p\times \zet_{2}\times\zet_{2}$ for which the quotient group $\Gamma/\Gamma_0$ has at least three involutions.  Lemma~\ref{p22}  implies that there exists a  MR$_{\Gamma_0}(p, 4)$  applying now Lemma~\ref{lemgl} we obtain a  MRS$_{\Gamma}(p, b;abc/4p)$. To  construct one of  the rectangles from MRS$_{\Gamma}(a, b;c)$ we take  $ab/4p$ of  MRS$_{\Gamma}(p, b;ac/4p)$ rectangles and ''glue'' them  into a  rectangle of size $a \times b$.

If  $c$ is odd, then $b =2^{\alpha_1+\ldots+\alpha_k}$ and moreover there exists a subgroup $\Gamma_0$ of $\Gamma$ such that $\Gamma_0\cong \zet_{p}\times  \zet_{2^{\alpha_1}}\times\zet_{2^{\alpha_2}}\times \ldots \times \zet_{2^{\alpha_k}}$. Observe that $|\Gamma/\Gamma_0|=ac/p$ is odd thus $\Gamma/\Gamma_0\in\gr$.   Lemma~\ref{p22}  implies that there exists a  MR$_{\Gamma_0}(p, b)$ and therefore using similar arguments as above we are done.~\qed\\

Moreover we are able to show the following.
\begin{lem}\label{p22_} 
 If $\Gamma\cong\zet_{2k+1}\times (\zet_{2})^{m}$  then  a  $\Gamma$-magic rectangle  set MRS$_{\Gamma}(2k+1, 2;2^{m-1})$ does not exist.
\end{lem}
\textit{Proof.} We can assume that $\Gamma\cong\zet_{2k+1}\times \Delta$ for $\Delta=\{0,\iota_1,\iota_2,\ldots,\iota_{2^{m}-1}\}$, where $\iota_j$ is an involution. Suppose that there exists  a  MRS$_{\Gamma}(2k+1, 2;2^{m-1})$.\\

 Note that  the entry in the $i$-th row and $j$-th column of the $s$-th rectangle in the $\Gamma$-magic rectangle set MRS$_{\Gamma}(2k+1, 2; 2^{m-1})$ has two coordinates $(i^s,j^s)$ for $i^s\in\zet_{2k+1}$ and $j^s \in \Delta$, moreover  all row sums in
every rectangle equal to a constant $(\omega_1,\omega_2)\in \zet_{2k+1}\times\Delta$  and all column sums in every rectangle equal to a
constant $(\delta_1, \delta_2)\in \zet_{2k+1}\times\Delta$. Note that every element from $\Delta$ is used odd times, therefore $\omega_2\neq0$ and without loss of generality we can assume that $\omega_2=\iota_1$. Hence for each row from the group $\Delta$ we are using one of the pair $\{0,\iota_1\}$, $\{\iota_j,\iota_l\}$ (for $\iota_j+\iota_l=\iota_1$) and each pair is used $2k+1$ times.
 One can easily see that it is impossible to obtain a constant column sum $\delta_2\in \Delta$.~\qed\\

We obtain the sufficient and necessary conditions for existence a $\Gamma$-magic rectangle MR$_{\Gamma}(a,b)$.
\begin{thm}\label{main3} Let   $a,b>1$.
A  $\Gamma$-magic rectangle  MR$_{\Gamma}(a, b)$ exists if and only if $a$ and $b$ are both even or $\Gamma\in\gr$.
\end{thm}
\textit{Proof.} Recall that MR$_{\Gamma}(a, b)$=MRS$_{\Gamma}(a, b;1)$. If $\{a,b\}= \{2,2m+1\}$ then a MR$_{\Gamma}(a, b)$ does not exist by Observation~\ref{odd} since any group $\Gamma$ of order $ab$ has exactly one involution. Therefore we are done by Theorems~\ref{main2} and \ref{a22}.~\qed\\

\section{Concluding remarks}\label{SectionOurConcludingRemarks}

Observe that by Observation~\ref{odd} and Theorems~\ref{main2} and \ref{a22} the   missing cases in the full characterization of $\Gamma$-magic rectangle sets  MRS$_{\Gamma}(a, b;c)$ are:
\begin{enumerate}
	\item $\{a,b\}=\{2k+1,2\}$, $c$ even  and $\Gamma$ has more than one involution;
	\item $\{a,b\}=\{2k+1,4\}$, $c\equiv 2 \pmod 4$ and $\Gamma\cong \Delta \times \zet_{2}\times\zet_{2}\times  \zet_{2}$ for  a group $\Delta$ of odd order $(2k+1)c/2$;
	\item $\{a,b\}=\{2k+1,2^\alpha\}$, $c= 2^{\beta}(2l+1)$ and $\Gamma\cong \Delta \times \zet_{2}\times\zet_{2^{\alpha+\beta-1}}$ for  a group $\Delta$ of odd order $(2k+1)(2l+1)$.
	
\end{enumerate}
We finish this section with a conjecture and some open problems.

\begin{conj} Let $a,b> 1$.
 A  $\Gamma$-magic rectangle set MRS$_{\Gamma}(a, b;c)$ exists if and only if $a$ and $b$ are both even or $\Gamma\in\gr$ and $\{a,b\}\neq\{2k+1,2\}$.
\end{conj}
Since in general row sums and column sums for  a $\Gamma$-magic rectangle set MRS$_{\Gamma}(a, b;c)$ are not unique  for the further research it can be also interesting to determine all possible row sums and column sums.

\begin{prob}
For  a $\Gamma$-magic rectangle set MRS$_{\Gamma}(a, b;c)$  determine all possible row sums and column sums.
\end{prob}
Moreover because for a $\Gamma$-magic rectangle set MRS$_{\Gamma}(a, b;c)$  we can define the order of the operation (summation) in a natural way (from left to right, and from up to down) hence we post also the following open problem.
\begin{prob}
Give necessary and sufficient conditions 
for  existence of a $\Gamma$-magic rectangle set MRS$_{\Gamma}(a, b;c)$, where $\Gamma$ is non Abelian group of oder $abc$.
\end{prob}


\bibliographystyle{plain}

\end{document}